\documentclass[11pt,a4papar]{article}
\usepackage{amsthm}
\usepackage{amssymb}
\usepackage{amscd}
\usepackage{amsmath}
\usepackage[all]{xy}
\usepackage{comment}
\usepackage{algorithmic,algorithm}
\usepackage{here}
\usepackage{graphicx}
\usepackage{color} %
\usepackage{enumerate} %
\usepackage{listliketab}
\usepackage{multirow}
\usepackage{cases}


\makeatletter
  
  \@addtoreset{algorithm}{subsection}
\makeatother
\newtheorem{thm}{Theorem}[section]
\newtheorem{crl}[thm]{Corollary}

\newtheorem{lmm}[thm]{Lemma}

\newtheorem{prp}[thm]{Proposition}

\theoremstyle{definition}
\newtheorem{dfn}[thm]{Definition}
\newtheorem{str}[thm]{Strategy}
\newtheorem{prob}[thm]{Problem}
\theoremstyle{remark}
\newtheorem*{rem}{Remark}%%% * means no numbering

\numberwithin{equation}{section}

\def\deg{{\rm deg}}

\begin{document}

%Title
\title{Counting isomorphism classes of superspecial curves}
\author{Momonari Kudo}
% \thanks{Kobe City College of Technology.}
% \thanks{Institute of Mathematics for Industry, Kyushu University. E-mail: \texttt{m-kudo@math.kyushu-u.ac.jp}}}
\providecommand{\keywords}[1]{\textbf{\textit{Key words---}} #1}
\maketitle

%\tableofcontents      %optional
\begin{abstract}      %optional
A superspecial curve is a (non-singular) curve over a field of positive characteristic whose Jacobian variety is isomorphic to a product of supersingular elliptic curves over the algebraic closure.
It is known that for given genus and characteristic, there exist only finitely many superspecial curves, up to isomorphism over an algebraically closed field.
In this article, we give a brief survey on results of counting isomorphism classes of superspecial curves.
In particular, this article summarizes some recent results in the case of genera four and five, obtained by the author and S.\ Harashita.
We also survey results obtained in a joint work with Harashita and E.\ W.\ Howe, on the enumeration of superspecial curves in a certain class of non-hyperelliptic curves of genus four.
\end{abstract}

% \footnote[0]{
% 2010 Mathematical Subject Classification: 14H45, 14G15, 14G05, 14Q05
% }
\footnote[0]{
This research was supported by JSPS Grant-in-Aid for Young Scientists 20K14301.
}

%=====================
\section{Introduction}
%=====================

Throughout this article, by a curve we mean a non-singular projective variety of dimension one.
A curve of genus $g$ over a field $K$ of characteristic $p>0$ is said to be {\it superspecial} if $\mathrm{Jac}(C) \cong E^g$ (over the algebraic closure $\overline{K}$) for a supersingular elliptic curve $E$, where $\mathrm{Jac}(C)$ denotes the Jacobian variety of $C$.
Note that this definition is well-defined by the following fact of Deligne, Ogus and Shioda (cf.\ \cite[Theorem 3.5]{Shioda} or \cite[Section 1.6, p.\ 13]{LO}):
If $g \geq 2$, for any supersingular elliptic curves $E_i$ for $1 \leq i \leq 2g$, we have $E_1 \times \cdots \times E_g \cong E_{g+1} \times \cdots \times E_{2g}$.

For a pair $(g,p)$, we denote by $\Lambda_{g,p}$ the set of $\overline{\mathbb{F}_{p}}$-isomorphism classes of superspecial curves of genus $g$ over finite fields of characteristic $p>0$.
The cardinality $\# \Lambda_{g,p}$ is at most finite (zero is possible) by, e.g., a general fact that given an abelian variety $A$, there exist only finitely many (irreducible) curves $D$ such that $\mathrm{Jac}(D) \cong A$, see \cite[Corollary 1.2]{NN}.
Ekedahl proved in \cite[Theorem 1.1]{Ekedahl}, which we will recall in Theorem \ref{thm:Ekedahl} in this article, that if there exists a superspeical curve $C$ of genus $g$ over $\overline{\mathbb{F}_p}$, then we have $2g \leq p^2 - p$, and $2g \leq p - 1$ if $C$ is hyperelliptic and $(g, p) \neq (1, 2)$.

The main problem which we consider in this article is:

\begin{prob}
Determine the number $\# \Lambda_{g,p}$ of $\overline{\mathbb{F}_{p}}$-isomorphism classes of superspecial curves of genus $g$ over finite fields of characteristic $p>0$.
Moreover, find complete representatives of the isomorphism classes.
\end{prob}

This article is a survey on results of counting the number of isomorphism classes of superspecial curves.
For $g=1$, Deuring~\cite{Deuring} proved that $\# \Lambda_{1,p}$ is equal to the class number of a quaternion algebra.
Also in the case of $g =2,3$, it follows from a general result by Ibukiyama-Katsura-Oort \cite[Theorem 2.10]{IKO} on superspecial principally polarized abelian varieties that $\# \Lambda_{g,p}$ is determined by computing the class numbers of quaternion hermitian lattices.
These class numbers were explicitly computed by Eichler~\cite{Eichler} for $g=1$, Hashimoto-Ibukiyama~\cite{HI} for $g=2$, and Hashimoto~\cite{H} for $g=3$ (cf.\ Igusa also computed the class number for $g=1$ by directly counting the $\overline{\mathbb{F}_p}$-isomorphism classes of supersingular elliptic curves). 
We review these results for $g \leq 3$ in Section \ref{sec:genus123}.

On the other hand, the problem for $g \geq 4$ has not been solved in all primes, but in recent years, the author and Harashita developed several algorithms to count genus-$4$ or genus-$5$ superspecial curves~\cite{KH17}, \cite{KH18}, \cite{KH20}, \cite{KH20b}.
Sections \ref{sec:genus4} and \ref{sec:genus5} describe our results (for small primes) obtained by these algorithms in the case of genus $g=4,5$.
In Section \ref{sec:app}, we also describe our most recent results obtained by a joint work with Harashita and Howe~\cite{KHH}, where we presented algorithms to count (or find) superspecial curves among certain $2$-dimensional families of genus-$4$ non-hyperelliptic curves. 

\begin{rem}
As stated above, this article mainly focuses on the enumeration of superspecial curves of given genus, up to isomorphism over {\it an algebraically closed field}, i.e., counting $\overline{\mathbb{F}_p}$-isomorphism classes of such curves over $\overline{\mathbb{F}_p}$.
On the other hand, we can also consider the enumeration up to isomorphism over {\it finite fields}, i.e., counting the number of $K$-isomorphism classes of superspecial curves over a finite field $K$.
Note that it suffices for this to consider the case of $K=\mathbb{F}_p$ or $\mathbb{F}_{p^2}$ since the number of $\mathbb{F}_{p^a}$-isomorphism classes of superspecial curves over $\mathbb{F}_{p^a}$ depends on the parity of $a$ (cf.\ \cite[Proposition 2.3.1]{KH20}).
Main results in \cite{KH17}, \cite{KH18}, \cite{KH20}, \cite{KH20b} that will be stated in Sections \ref{sec:genus4} and \ref{sec:genus5} include those on the enumeration not only over $\overline{\mathbb{F}_{p}}$, but also over $\mathbb{F}_{p}$ or $\mathbb{F}_{p^2}$.
\end{rem}

%==========================================================================================
\section{The number of superspecial curves of genus one, two and three}\label{sec:genus123}
%==========================================================================================

Let $p$ be a rational prime.
At first, we recall Ekedahl's results in \cite{Ekedahl} including the field of definition of superspecial curves.
The main theorem of \cite{Ekedahl} (Theorem \ref{thm:Ekedahl} below) gives bounds on the existence of superspecial curves:

\begin{thm}[\cite{Ekedahl}, Theorem 1.1]\label{thm:Ekedahl}
If there exists a superspecial curve $C$ of genus $g$ in characteristic $p$, then we have the following:
\begin{enumerate}
\item $2g \leq p^2 - p$, and
\item $2g \leq p - 1$ if $C$ is hyperelliptic and $(g, p) \neq (1, 2)$.
\end{enumerate}
\end{thm}

Ekedahl also showed in the proof of \cite[Theorem 1.1]{Ekedahl} that any superspecial curve over an algebraically closed field descends to a maximal or minimal curve over $\mathbb{F}_{p^2}$, where a curve $C$ of genus $g$ over $\mathbb{F}_q$ is called {\it maximal} (resp. {\it minimal}) if the number of $\mathbb{F}_q$-rational points on $C$ attains the Hasse-Weil upper bound $q + 1 + 2 g \sqrt{q}$ (resp. the Hasse-Weil lower bound $q + 1 - 2 g \sqrt{q}$).
Conversely, it is known that any maximal or minimal curve over $\mathbb{F}_{p^2}$ is superspecial.
Thus, for determining $\# \Lambda_{g,p}$, it suffices to count $\overline{\mathbb{F}_{p}}$-isomorphism classes of superspecial curves of genus $g$ over $\mathbb{F}_{p^2}$.
One more important fact showed in the proof of \cite[Theorem 1.1]{Ekedahl} is that the existence of a superspecial curve over the prime field $\mathbb{F}_{p}$ implies that of maximal and minimal curves over $\mathbb{F}_{p^2}$.

%=============================================================
% \subsection{Genus $g=1$: Elliptic curves}\label{subsec:genus1}
%=============================================================
In the following, we review results on the computation of $\# \Lambda_{g,p}$ for $g \leq 3$.
Let $h_p$ (resp.\ $t_p$) denote the class (resp.\ type) number of the quaternion algebra $B_{p,\infty}$ over $\mathbb{Q}$ ramified exactly at $\{ p, \infty \}$.
Deuring~\cite{Deuring} showed that the computation of $\# \Lambda_{1,p}$ is reduced into that of the class number $h_p$.

\begin{thm}[\cite{Deuring}]\label{thm:deuring}
We have the following:
\begin{enumerate}
\item Every supersingular elliptic curve over $\overline{\mathbb{F}_p}$ has a model over $\mathbb{F}_{p^2}$, and $\# \Lambda_{1,p}=h_p$.
\item The number of elements in $\Lambda_{1,p}$ which have models defined over $\mathbb{F}_p$ is $2 t_p - h_p$.
\end{enumerate}
\end{thm}

Using results on computing $h_p$ and $t_p$ by Eichler~\cite{Eichler}, we have the following:

\begin{thm}[\cite{Deuring}, \cite{Eichler}]\label{thm:eichler}
The number $\# \Lambda_{1,p}$ of $\overline{\mathbb{F}_p}$-isomorphism classes of supersingular elliptic curves is equal to
\[
\frac{p-1}{12} + \frac{1 - \left( \frac{-1}{p} \right) }{4} + \frac{1 - \left( \frac{-3}{p} \right) }{3}
\] 
if $p >3$, and one if $p=2$ or $3$.
\end{thm}

Igusa~\cite{Igusa} also proved the same result as in Theorem \ref{thm:eichler} by directly computing the number of supersingular $j$-invariants from the Legendre form $y^2 = x(x-1)(x-\lambda)$ of an elliptic curve.
We also refer to \cite[Proposition 4.4]{XYY16} for results on the number of $\mathbb{F}_q$-isomorphism classes of supersingular elliptic curves over $\mathbb{F}_q$.

%============================
% \subsection{Genus two and three}
%============================

For $g=2$ and $3$, determining $\# \Lambda_{g,p}$ is reduced into counting superspecial principally polarized abelian varieties (PPAV's for short) by the fact that any PPAV is the Jacobian variety of a (possibly reducible) curve, see the main theorem of \cite{OU}.
Here we recall a general result by Ibukiyama-Katsura-Oort \cite[Theorem 2.10]{IKO} on the number of superspecial principally polarized abelian varieties:

\begin{thm}[\cite{IKO}, Theorem 2.10]\label{thm:IKO}
Let $E$ be a supersingular elliptic curve over $\overline{\mathbb{F}_p}$.
For $g\geq 2$, the number of principal polarizations on $E^g$ up to automorphisms of $E^g$ is equal to the class number $H_g = H_g(p, 1)$ of the principal genus of the quaternion hermitian space $(B_{p,\infty})^g$.
\end{thm}

% the dimension of the moduli space of principally polarized abelian varieties of dimension $g$ is equal to that of the moduli space of curves of genus $g$.
Counting superspecial curves of genus $g =2,3$ is done by removing the contribution of reducible curves.
As is noted in \cite[p.\ 145]{IKO}, the number of supersingular abelian surfaces with reducible principal polarization is equal to the number of pairs $(E_i,E_j)$ of supersingular elliptic curves $E_i$ and $E_j$ with $i \leq j$, and thus we have the following result for $g=2$:

\begin{thm}[\cite{IKO}, \cite{HI}]\label{thm:genus2}
The number $\# \Lambda_{2,p}$ of $\overline{\mathbb{F}_p}$-isomorphism classes of superspecial curves of genus two is equal to $H_2 - H_1(H_1+1)/2$.
Using the computational result of $H_2$ by Hashimoto-Ibukiyama~\cite{HI} together with Theorem \ref{thm:eichler}, we have that
\[
\# \Lambda_{2,p} =
\begin{cases}
0 & (p=2,3) \\
1 & (p=5) \\
 \frac{p^3 + 24 p^2 + 141 p - 166}{2880} - \frac{1 - \left( \frac{-1}{p} \right) }{32} + \frac{1 - \left( \frac{-2}{p} \right) }{8} + \frac{1 - \left( \frac{-3}{p} \right) }{18} + \epsilon & (p \geq 7),
 \end{cases}
\] 
where $\epsilon = 4/5$ if $p \equiv 4 \pmod{5}$, and zero otherwise.
\end{thm}

The number of $\overline{\mathbb{F}_p}$-isomorphism classes of superspecial curves $C$ of genus $2$ such that $C$ has a model over $\mathbb{F}_p$ is also computable, see \cite[Section 1]{IK}.

Similarly to the case of $g=2$, we can compute the value of $\# \Lambda_{3,p}$, see \cite[Theorem 3.10 (d)]{Brock} for an explicit formula.

% \begin{thm}[\cite{IKO}, \cite{H}]\label{thm:genus3}
% The number $\# \Lambda_{3,p}$ of $\overline{\mathbb{F}_p}$-isomorphism classes of superspecial curves of genus three is equal to $H_3 - H_1 H_2 - H_1 (H_1 + 1) (H_2 + 2) /6$.
% Using the computational result of $H_3$ by Hashimoto~\cite{H}
% \[
% \frac{p^3 + 24 p^2 + 141 p - 166}{2880} - \frac{1 - \left( \frac{-1}{p} \right) }{32} + \frac{1 - \left( \frac{-2}{p} \right) }{8} + \frac{1 - \left( \frac{-3}{p} \right) }{18} 
% \] 
% if $p >5$, one if $p=5$, and zero if $p=2,3$.
% \end{thm}

%===================================================================
\section{Case of genus four}\label{sec:genus4}
%===================================================================

Different from the case of $g \leq 3$, for $g \geq 4$ the dimension of the moduli space of curves of genus $g$ is strictly less than that of the moduli space of PPAV's of dimension $g$. 
This means that the theory of abelian varieties is not so effective for our purpose for $g \geq 4$.
For this reason, the enumeration of superspecial curves of genus $4$ has not been completed yet for every $p$, whereas some results for small and concrete $p$ are known.
In this section (resp.\ the next section), we survey results on the enumeration of superspecial curves of genus $4$ (resp.\ $5$).
In particular, this section summarizes results in \cite{KH17}, \cite{KH18} and \cite{KH20}, where the authors proposed computational approaches to enumerate superspecial curves of genus $4$.

Let $C$ be a curve of genus $4$.
We recall that $C$ is either of the following two types (cf.\ \cite[Chap.\ IV, Example 5.2.2]{Har}):
\begin{enumerate}
\item {\it Hyperelliptic}. The normalization of the plane curve $y^2=f(x)$, where $f(x)$ is a separable polynomial of degree $9$ or $10$.
\item {\it Canonical}. A complete intersection of quadratic and cubic hypersurfaces in $\mathbb{P}^3$.
\end{enumerate}

Recall from the paragraph just after Theorem \ref{thm:Ekedahl} that, for counting superspecial curves in characteristic $p>0$, it suffices to count $\overline{\mathbb{F}_p}$-isomorphism classes of superspecial curves over $\mathbb{F}_{p^2}$.

First, we consider the case where $C$ is a non-hyperelliptic curve over a finite field $K=\mathbb{F}_q$ with $q = p$ or $p^2$ for $p \geq 5$, and give a summary of results in \cite{KH17} and \cite{KH20}.
As a canonical curve, $C$ is defined in the $3$-projective space $\mathbb{P}^3=\mathrm{Proj}(\overline{K}[x,y,z,w])$ by an irreducible quadratic form $Q$ and an irreducible cubic form $P$ in $\overline{K}[x,y,z,w]$, see \cite[Chapter IV, Example 5.2.2]{Har}.
As showed in \cite[Section 2.1]{KH17}, we may assume that any coefficient of $Q$ and $P$ belongs to $K$.
By the classification theory of quadratic forms over finite fields, we can transform $Q$ into either of {\bf (N1)} $2xw + 2yz$, {\bf (N2)} $2xw + y^2 - \epsilon z^2$ for $\epsilon\in K^\times \smallsetminus (K^\times)^2$ and {\bf (Dege)} $2yw + z^2$ (cf.\ \cite[Remark 2.1.1]{KH17}).

Here we recall a criterion on the superspecialty of $C=V(Q,P)$ in Proposition \ref{cor:HW} below.
The proof is done by computing the Hasse-Witt matrix of $C$, which represents the Frobenius on the first cohomology group $H^1 (C, \mathcal{O}_C)$.
Each entry of the Hasse-Witt matrix of $C$ is one of the $16$ coefficients in $(QP)^{p-1}$ given in Proposition \ref{cor:HW}.

\begin{prp}[\cite{KH17}, Corollary 3.1.6]\label{cor:HW}
With notation as above, $C$ is superspecial if and only if the coefficients of $x^{pi-i'}y^{pj-j'}z^{pk-k'}w^{p \ell-\ell'}$
in $(Q P )^{p-1}$ are equal to $0$
for all positive integers $i,j,k,\ell,i',j',k',\ell'$ with $i+j+k+\ell=i'+j'+k'+\ell'=5$.
\end{prp}

Based on Proposition \ref{cor:HW}, we have a computational strategy to enumerate superspecial non-hyperelliptic curves of genus $4$ over $K$:
\begin{str}\label{strategy1}
\begin{enumerate}
\item For each of the three types ({\bf (N1)}, {\bf (N2)} and {\bf (Dege)}) of $Q$, collect superspecial curves $V(Q,P)$ as follows: 
\begin{enumerate}
\item Collect cubic forms $P \in K[x,y,z,w]$ such that $\mathrm{HW}\text{-matrix} = 0$, i.e., the $16$ coefficients in $(QP)^{p-1}$ given in Proposition \ref{cor:HW} are all zero. 
\item For each $P$ collected in (a), test whether $V(Q,P)$ is non-singular or not.
\end{enumerate}
\item For the superspecial curves $V(Q,P)$ collected in Step 1, compute their isomorphism classes.
\end{enumerate}
\end{str}

Both of Step 1 (a), (b) and Step 2 are done with Gr\"{o}bner basis computation.
We here focus on Step 1 (a) and Step 2, and give their brief descriptions.
Note that Step 1 (b) is done by a general method for the non-singularity test, see e.g., \cite[Section 3.2]{KH17}.

Step 1 computes the solutions of multivariate systems $\mathrm{HW}\text{-matrix} = 0$ with respect to unknown coefficients in $P$.
Naively, $P$ has $20$ unknowns, but in fact, the dimension of the moduli space of non-hyperelliptic curves of genus $4$ is $9$.
Thus it requires to reduce the number of unknowns as much as possible, since the number deeply affects the computational cost of solving multivariate systems.
The author and Harashita~\cite[Section 4]{KH17}, \cite[Section 3]{KH20} reduced the number by considering the action of elements in the orthogonal similitude group $\tilde{\mathrm{O}}_\varphi(K)$ to the cubic form $P$, where $\varphi$ is the symmetric matrix associated to $Q$.
Note that they realized elements in $\tilde{\mathrm{O}}_\varphi(K)$ by computing the Bruhat decomposition of $\tilde{\mathrm{O}}_\varphi(K)$, see \cite[Section 3]{KH17} for more details.

In Lemmas \ref{NewReductionLemmaN1} -- \ref{ReductionLemmaDegenerate} below, we collect the reduced form of $P$ for each of the three types of $Q$.

\begin{lmm}[\cite{KH20}, Lemma 3.4.1]\label{NewReductionLemmaN1}
Let $Q = 2 x w + 2 yz$, and $\varphi$ the symmetric matrix associated to $Q$.
An element of $\tilde{\mathrm{O}}_\varphi(K)$ transforms $P$ into the following form:
\begin{eqnarray}\label{eq:ReductionLemmaN1}
\begin{split}
& (y + b_1 z)x^2  + b_2 xz^2 + a_1 y^3 + a_2 y^2z + a_3 yz^2 +  a_4 z^3 \\
&  + (a_5 y^2 + a_6 yz + a_7 z^2)w + (a_8y + a_9z)w^2 + a_{10}w^3
\end{split}
\end{eqnarray}
for $a_i \in K$
and for $b_1 \in \{ 0 \} \cup K^\times/(K^\times)^2$ and $b_2 \in \{ 0, 1 \}$.
\end{lmm}

\begin{lmm}[\cite{KH20}, Lemma 3.5.1]\label{ReductionLemmaN2}
Let $Q = 2 x w + y^2 - \epsilon z^2$ for $\epsilon \in K^{\times}$ with $\epsilon \notin (K^{\times})^2$.
An element of $\tilde{\mathrm{O}}_\varphi(K)$ transforms $P$ into the following form:
\begin{eqnarray}\label{eq:ReductionLemmaN2}
\begin{split}
& (a_1 y + a_2 z)x^2  + a_3 (y^2 - \epsilon z^2)x
+ b_1 y (y^2 - \epsilon z^2) + a_4 y ( y^2 + 3 \epsilon z^2) \\
&  + a_5 z (3 y^2+\epsilon z^2) + (a_6 y^2 + a_7 y z + b_2 z^2)w + (a_8 y + a_9 z)w^2 + a_{10}w^3
\end{split}
\end{eqnarray}
for some $a_i \in K$ with $(a_1,a_2) \ne (0,0)$ and for $b_1, b_2 \in \{ 0, 1 \}$.
\end{lmm}

\begin{lmm}[\cite{KH20}, Lemma 3.6.1]\label{ReductionLemmaDegenerate}
Let $Q = 2 y w + z^2$, and $\varphi$ the symmetric matrix associated to $Q$.
An element of $\tilde{\mathrm{O}}_\varphi(K)$ transforms $P$ into the following form \eqref{eq1:ReductionLemmaDegenerate} if $\# K > 5$, and into either of the following forms \eqref{eq1:ReductionLemmaDegenerate} and \eqref{eq2:ReductionLemmaDegenerate} if $\# K = 5$:
\begin{eqnarray}\label{eq1:ReductionLemmaDegenerate}
\begin{split}
& a_0 x^3 + (a_1 y^2 + a_2 z^2 + a_3 w^2 + a_4 y z + a_5 z w)x\\
& + a_6 y^3 + a_7 z^3 + a_8 w^3 + a_9 y z^2 + b_1 z^2 w + b_2 z w^2,
\end{split}
\end{eqnarray}
for some $a_i\in K$ with $a_0, a_6 \in K^\times$ and for $b_1, b_2 \in \{0,1\}$, where the leading coefficient of $R := a_1 y^2 + a_2 z^2 + a_3 w^2 + a_4 y z + a_5 z w$ is $1$ or $R=0$;
\begin{equation}\label{eq2:ReductionLemmaDegenerate}
x^3 + (a_1 y^2 + a_2 z^2 + a_3 w^2 + a_4 y z + b_1 z w)x + y^2 z + z w^2
\end{equation}
for $a_i \in K=\mathbb{F}_5$ and $b_1 \in \{0, 1\}$.
\end{lmm}

Thus Step 1 (a) of Strategy \ref{strategy1} is done by:
For each of the three types of $Q$, collect $P$ in the corresponding reduced form given in Lemma \ref{NewReductionLemmaN1} for {\bf (N1)}, Lemma \ref{ReductionLemmaN2} for {\bf (N2)} and Lemma \ref{ReductionLemmaDegenerate} for {\bf (Dege)} such that $\mathrm{HW}\text{-matrix} = 0$.

We next consider Step 2 of Strategy \ref{strategy1}, i.e., how to decide whether two non-hyperelliptic curves of genus $4$ are isomorphic or not over $k$, where $k = K$ or $\overline{K}$.
Let $C_1 = V (Q_1, P_1)$ and $C_2 = V(Q_2, P_2)$ be non-hyperelliptic curves of genus $4$ over $k$.
If there exists an isomorphism over $k$ from $C_1$ to $C_2$, the quadratic forms $Q_1$ and $Q_2$ are equivalent over $k$.
Hence it suffices to consider the case of $Q_1=Q_2$, say $Q$.
Let $\varphi$ be the symmetric matrix associated to $Q$.
As is described in \cite[Section 6.1]{KH17}, the two curves $C_1$ and $C_2$ are $k$-isomorphic if and only if there exist $g \in \tilde{\mathrm{O}}_\varphi(k) $ and $\lambda \in k^{\times}$ such that $g \cdot P_1 \equiv \lambda P_2 \mbox{ mod } Q$.
Using the Bruhat decomposition of $\tilde{\mathrm{O}}_\varphi(k) $ given in \cite[Section 3]{KH17}, we reduce the (non-)existence of such $g$ and $\lambda$ into that of solutions over $k$ of multivariate systems, see \cite[Section 4.2]{KH20} for more details including concrete algorithms to test the (non)-existence of such solutions.

In \cite[Section 5]{KH17} and \cite[Section 4]{KH20}, the authors wrote down explicit algorithms for Strategy \ref{strategy1}, and implemented them over Magma~\cite{Magma}.
To implement Step 1 (a), they adopted the hybrid method~\cite{BFP}, which combines the Gr\"{o}bner basis computation with the brute force on some unknown coefficients in $P$.
In Theorem \ref{MainTheorem} below, we collect main results in \cite{KH17} and \cite{KH20} obtained by executing proposed algorithms over Magma.

\begin{thm}\label{MainTheorem}
\begin{enumerate}
	\item (\cite[Theorem A]{KH17}) Any superspecial curve of genus $4$ over $\mathbb{F}_{5^2}$ is $\mathbb{F}_{5^2}$-isomorphic to $2 y w + z^2 = x^3 + a_1 y^3 + a_2 w^3 + a_3 z w^2 = 0$ in $\mathbb{P}^3$, where $a_1,a_2 \in \mathbb{F}_{5^2}^{\times}$ and $a_3 \in \mathbb{F}_{5^2}$.
	\item (\cite[Corollary 5.1.1]{KH17}) All superspecial curves of genus $4$ in characteristic $5$ are isomorphic to each other over an algebraically closed field.
	\item (\cite[Theorem A]{KH20} and \cite[Example 6.2.4]{KH17}) There exist exactly seven (resp.\ $21$) superspecial curves of genus $4$ over $\mathbb{F}_{5}$ (resp.\ $\mathbb{F}_{25}$) up to isomorphism over $\mathbb{F}_{5}$ (resp.\ $\mathbb{F}_{25}$).
	\item (\cite[Theorem B]{KH17}) There is no superspecial curve of genus $4$ in characteristic $7$.
	\item (\cite[Theorem B]{KH20}) There exist exactly $30$ (resp.\ nine) non-hyperelliptic superspecial curves of genus $4$ over $\mathbb{F}_{11}$ up to isomorphism over $\mathbb{F}_{11}$ $($resp.\ $\overline{\mathbb{F}_{11}})$.
\end{enumerate}
\end{thm}

% \begin{table}[htb]
%   \begin{tabular}{|c||c|c|c|c|c|} \hline
%   $K$                       & $2$ & $2^2$ & $3$ & $3^2$ &  \\ \hline
%   $\mathbb{F}_q$            & $0$ & $0$   & $0$ & $0$   &  \\ \hline
%   $\overline{\mathbb{F}_q}$ & $0$ & $0$   & $0$ & $0$   &
%   \end{tabular}
% \end{table}

Second, we consider the case where $C$ is hyperelliptic, and give a summary of results in \cite{KH18}.
In general, a hyperelliptic curve $H$ over $K$ is realized as the desingularization of the homogenization of $y^2 = f(x)$, where $f(x)$ is a polynomial over $K$ with non-zero discriminant.
In \cite[Section 3.2]{KH18}, the authors gave a reduction of a defining equation of $H$ so that the set of all the ramification points of the reduced model is defined over $K$:

\begin{lmm}[\cite{KH18}, Lemma 2]\label{ReductionHyper}
Assume that $p$ and $2g+2$ are coprime.
Let $\epsilon \in K^\times \smallsetminus (K^\times)^2$.
Any hyperelliptic curve $H$ of genus $g$ over $K$ is the desingularization of
the homogenization of
\[
c y^2 = x^{2g+2} + b x^{2g} + a_{2g-1}x^{2g-1} + \cdots + a_1 x + a_0
\]
for $a_i \in K$ for $i=0,1,\ldots, 2g-1$ where $b= 0, 1,\epsilon$ and $c=1,\epsilon$.
\end{lmm}

We also recall a criterion on the superspecialty of the hyperelliptic curve $H$ of genus $g$ in Corollary \ref{cor:HWh} below.
This criterion comes from a well-known explicit formula for the Cartier-Manin matrix of $H$ (cf.\ \cite{Yui}), which represents the Cartier operator on the space $H^0 (H, \Omega_H^1)$ of regular differential forms on $H$. 

\begin{crl}[\cite{KH18}, Corollary 1]\label{cor:HWh}
Let $H$ be a hyperelliptic curve $y^2 = f(x)$ of genus $g$ over $K$, where $\deg (f) = 2 g + 1$ or $2 g+2$.
Then $H$ is superspecial if and only if the coefficients of $x^{p i - j}$ in $f^{(p-1)/2}$ are equal to $0$ for all pairs of integers with $1 \leq i , j \leq g$.
\end{crl}

Here we also describe a method given in \cite[Section 3.3]{KH18} to test whether two hyperelliptic curves $C_1 : c_1 y^2 = f_1(x)$ and $C_2 : c_2 y^2 = f_2(x)$ with $c_1,c_2 \in K^{\times}$ of genus $g$ are isomorphic or not over $k=K$ or $\overline{K}$, where $f_i$ is a separable polynomial in $K[x]$ of degree $2 g + 2$ for each $1\leq i \leq 2$.
Let $F_i$ be the homogenization of $f_i$ with respect to an extra variable $z$ for each $1 \leq i \leq 2$.
Recall from \cite[Lemma 1]{KH18} that $C_1 \cong C_2$ over $k$ if and only if there exist $h \in \mathrm{GL}_2 ( k )$ and $\lambda \in k^{\times}$ such that $h \cdot F_1 = \lambda^2 F_2$.
Regarding entries of $h$ and $\lambda$ as variables, we reduce the (non-)existence of such $h$ and $\lambda$ into that of a solution over $k$ of a multivariate system, see \cite[Section 3.3]{KH18} for more details including concrete algorithms to test the (non)-existence of such a solution.

Combining Lemma \ref{ReductionHyper} and Corollary \ref{cor:HWh} with the isomorphism test described as above, we can construct a strategy similar to Strategy \ref{strategy1} for enumerating superspecial hyperelliptic curves of genus $g$ over $K$, see \cite[Section 3]{KH18} for concrete algorithms.
The authors of \cite{KH18} implemented the algorithms over Magma, and executed them for $g=4$ with $q=11,11^2,13,13^2,17,17^2,19$.
Theorem \ref{MainTheorem12} below collects main results in \cite{KH18}.

\begin{thm}\label{MainTheorem12}
\begin{enumerate}
	\item (\cite[Theorem 1]{KH18}) There is no superspecial hyperelliptic curve of genus $4$ in characteristic $11$ and $13$.
	\item (\cite[Theorem 2]{KH18}) There exist precisely $5$ $($resp.\ $25)$ superspecial hyperelliptic curves of genus $4$ over $\mathbb{F}_{17}$ $($resp.\ $\mathbb{F}_{17^2})$ up to isomorphism over $\mathbb{F}_{17}$ $($resp.\ $\mathbb{F}_{17^2})$.
Moreover, there exist precisely $2$ superspecial hyperelliptic curves of genus $4$ over the algebraic closure in characteristic $17$ up to isomorphism.
	\item (\cite[Theorem 3]{KH18}) There exist precisely $12$ superspecial hyperelliptic curves of genus $4$ over $\mathbb{F}_{19}$ up to isomorphism over $\mathbb{F}_{19}$.
Moreover, there exist precisely $2$ superspecial hyperelliptic curves of genus $4$ over $\mathbb{F}_{19}$ up to isomorphism over the algebraic closure.
\end{enumerate}
\end{thm}

Table \ref{tab:genus4} summarizes known values of the number $\# \Lambda_{4,p}$ of $\overline{\mathbb{F}_p}$-isomorphism classes of superspecial curves of genus $4$ in characteristic $p$.
The non-existence of non-hyperelliptic (resp.\ hyperelliptic) superspecial curves of genus $4$ for $p \leq 3$ (resp.\ $p \leq 7$) is deduced from Ekedahl's bounds given in Theorem \ref{thm:Ekedahl}.
The number written in bold type is determined by our theorems (Theorems \ref{MainTheorem} and \ref{MainTheorem12}) described in this section.
The notation `{\bf H}' and `{\bf C}' denote the hyperelliptic and canonical cases respectively.
The number written in each bracket is the number of $\overline{\mathbb{F}_p}$-isomorphism classes of  superspecial curves $C$ such that $C$ has a model over $\mathbb{F}_p$.

\begin{table}[htb]
\centering
    \caption{Known values of the number $\# \Lambda_{4,p}$ of $\overline{\mathbb{F}_p}$-isomorphism classes of superspecial curves of genus $4$ in characteristic $p$.
}\label{tab:genus4}
  \begin{tabular}{|c||c|c|c|c|c|c|c|c|c|} \hline
  $p$                    & $2$   & $3$   & $5$  & $7$   & $11$     &  $13$ & $17$  & $19$ & $\geq 23$\\ \hline
  \multirow{2}{*}{\bf H} & \multirow{2}{*}{0} & \multirow{2}{*}{0} & \multirow{2}{*}{0} & \multirow{2}{*}{0} & \multirow{2}{*}{\bf 0} & \multirow{2}{*}{\bf 0} & {\bf 2} &   ?  & \multirow{2}{*}{?} \\ 
                         &  & & & & & & ({\bf 2})& ({\bf 2})&  \\ \hline 
  \multirow{2}{*}{\bf C} & \multirow{2}{*}{0}   & \multirow{2}{*}{0} &  {\bf 1}  & \multirow{2}{*}{\bf 0}   & ? &  \multirow{2}{*}{?} & \multirow{2}{*}{?} & \multirow{2}{*}{?}    & \multirow{2}{*}{?}\\ 
                         &  & & ({\bf 1}) & &({\bf 9})     &      & & &  \\ \hline 
  \multirow{2}{*}{$\# \Lambda_{4,p}$} & \multirow{2}{*}{0}   & \multirow{2}{*}{0} &  {\bf 1}  & \multirow{2}{*}{\bf 0}   & ? &  \multirow{2}{*}{?} & \multirow{2}{*}{?} & \multirow{2}{*}{?}    & \multirow{2}{*}{?}\\ 
                         &  & & ({\bf 1}) & &({\bf 9})     &      & & &  \\ \hline 
  \end{tabular}
\end{table}

We close this section with open problems in the enumeration of superspecial curves of genus four:

\begin{prob}[Genus four]
Determine the number of $K$ or $\overline{K}$-isomorphism classes of superspecial curves of genus four over $K$ in the following cases:
\begin{enumerate}
\item Canonical case over $K = \mathbb{F}_{p}$ for $p \geq 13$ or over $K =\mathbb{F}_{p^2}$ for $p \geq 11$.
\item Hyperelliptic case over $K = \mathbb{F}_{p}$ for $p \geq 23$ or over $K =\mathbb{F}_{p^2}$ for $p \geq 19$.
\end{enumerate}
\end{prob}

%===================================================================
\section{Case of genus five}\label{sec:genus5}
%===================================================================

First, we recall that a curve of genus $5$ is either of the following three types:
\begin{enumerate}
\item {\it Hyperelliptic}. The normalization of the plane curve $y^2=f(x)$, where $f(x)$ is a separable polynomial of degree $11$ or $12$.
\item {\it Trigonal}. A curve $C$ such that there exists a morphism $C \rightarrow \mathbb{P}^1$ of degree $3$.
\item {\it Generic (canonical and non-trigonal)}. A complete intersection of three quadric in $\mathbb{P}^4$.
\end{enumerate}
In this paper, we say that a curve of genus $5$ is ``generic'' if it is canonical and non-trigonal.

As for the non-existence of superspecial curves of genus $5$, it follows from Ekedahl's bound given in Theorem \ref{thm:Ekedahl} that there is no superspecial curve of genus $5$ in characteristic $p=2$, $3$.
The non-existence holds also for $p=5$.
Indeed, by \cite[Lemma 2.2.1]{KH17}, if there were a superspecial curve of genus $5$ in characteristic $5$, then there would exist a maximal curve of genus $5$ over $\mathbb{F}_{5^2}$, which contradicts the fact due to Fuhrmann and Torres \cite{FT} that if there exists a maximal curve of genus $g$ over $\mathbb{F}_{p^2}$, then $4g \leq (p-1)^2$ or $2g=p^2-p$.
For $p \geq 7$, there is no result that shows the non-existence in the canonical case.
The problem of counting superspecial curves of genus $5$ is left for $p \geq 11$ in the hyperelliptic case, and for $p \geq 7$ in the canonical case.

This section briefly describes results in \cite{KH20b}, where the authors enumerated superspecial trigonal curves of genus $g=5$ over finite fields $\mathbb{F}_{p^a}$ for any $a$ if $p \leq 7$ and for odd $a$ if $p \leq 13$.
Recall from \cite[Proposition 2.3.1]{KH17} that it suffices to study the case of $a=1$, $2$ if $p \leq 7$ and the case of $a=1$ if $p \leq 13$.

Let $C$ be a trigonal curve of genus $5$ over a finite field $K=\mathbb{F}_q$, where $q = p$ or $p^2$ with $p \geq 5$.
It is shown in \cite[Section 2]{KH20b} that $C$ is the normalization of a plane quintic $C'=V(F) \subset \mathbb{P}^2$ with a unique singular point for some quintic form $F \in K[x,y,z]$.
The quintic forms defining our curves are divided into the following three types: {\bf (Split node case)} $F=xyz^3+f$, {\bf (Non-split node case)} $F=(x^2-\epsilon y^2)z^3+f$ with $\epsilon \in K \smallsetminus (K^\times)^2$, {\bf (Cusp case)} $F=x^2z^3+f$, where $f$ is the sum of monomial terms which can not be divided by $z^3$.
Specifically, we have the reduced forms of $F$ in Propositions \ref{ReductionSplitNode} -- \ref{ReductionCusp} below.
In the following, let $\zeta$ be a primitive element of $K^{\times}$, and $\epsilon$ an element of $K^\times \smallsetminus (K^\times)^2$.

\begin{prp}[\cite{KH20b}, Proposition 3.1.1]\label{ReductionSplitNode}
Any trigonal curve of genus $5$ over $K$ in {\bf (Split node case)} has a quintic model in $\mathbb{P}^2$ of the form \eqref{SplitNodeReducedEq1} or \eqref{SplitNodeReducedEq2}:
\begin{equation}
\begin{split}
F & = x y z^3 + (x^3 + b_1 y^3) z^2
+ (a_1 x^4 + a_2 x^3 y + a_3 x^2 y^2 + a_4 x y^3 + a_5 y^4) z \\
& \quad + a_6 x^5 + a_7 x^4 y + a_8 x^3 y^2 + a_9 x^2 y^3 + a_{10} x y^4 + a_{11} y^5, \label{SplitNodeReducedEq1}
\end{split}
\end{equation}
for $a_i \in K$, where $b_1 \in \{0, 1\}$ if $q \equiv 2 \bmod{3}$ and $b_1 \in \{0, 1, \zeta \}$ if $q \equiv 1 \bmod{3}$.
\begin{equation}
\begin{split}
F &= x y z^3 + (c_1 x^4 + c_2 x^3 y + a_3 x^2 y^2 + a_4 x y^3 + a_5 y^4) z\\
&\quad + a_6 x^5 + a_7 x^4 y + a_8 x^3 y^2 + a_9 x^2 y^3 + a_{10} x y^4 + a_{11} y^5. \label{SplitNodeReducedEq2}
\end{split}
\end{equation}
for $(c_1,c_2) = (0,0), (1,0), (0,1), (1,1), (1, \zeta)$ and for $a_i \in K$.
\end{prp}

\begin{prp}[\cite{KH20b}, Proposition 3.2.1]\label{ReductionNonSplitNode}
Any trigonal curve of genus $5$ over $K$ in {\bf (Non-split node case)} has a quintic model in $\mathbb{P}^2$ of the form \eqref{NonSplitNodeReducedEq1}, \eqref{NonSplitNodeReducedEq2} or \eqref{NonSplitNodeReducedEq3}:
\begin{equation}
\begin{split}
F &= (x^2-\epsilon y^2)z^3 + \{x(x^2+3\epsilon y^2) + b y(3x^2+\epsilon y^2)\} z^2\\
&\quad + (a_1 x^4 + a_2 x^3 y + a_3 x^2 y^2 + a_4 x y^3 + a_5 y^4) z\\
&\quad + a_6 x^5 + a_7 x^4 y + a_8 x^3 y^2 + a_9 x^2 y^3 + a_{10} x y^4 + a_{11} y^5, \label{NonSplitNodeReducedEq1}
\end{split}
\end{equation}
for $a_i\in K$, where $b=0$ if $q \not\equiv -1 \bmod{3}$ and otherwise $b$ has three possibilities determined by the condition that $(1,b)$ is parallel to $(1,0)A$ for a representative $A$
of $\tilde {\rm C}/\tilde {\rm C}^3$
{\rm (}for example $b=0,6,10$ if $q=11${\rm )}, where
\[
\tilde {\rm C} = \left\{\left.
\begin{pmatrix}r & \epsilon s\\ s & r\end{pmatrix}\right| (r,s) \in K^2,\ (r,s)\ne (0,0) \right\}.
\]
\begin{equation}
\begin{split}
F &= (x^2-\epsilon y^2)z^3 + (c x^4 + a_2 x^3 y + a_3 x^2 y^2 + a_4 x y^3 + a_5 y^4) z\\
&\quad + a_6 x^5 + a_7 x^4 y + a_8 x^3 y^2 + a_9 x^2 y^3 + a_{10} x y^4 + a_{11} y^5. \label{NonSplitNodeReducedEq2}
\end{split}
\end{equation}
for $c = 1, \zeta$ and for $a_i\in K$.
\begin{equation}
F = (x^2-\epsilon y^2)z^3 + a_6 x^5 + a_7 x^4 y + a_8 x^3 y^2 + a_9 x^2 y^3 + a_{10} x y^4 + a_{11} y^5. \label{NonSplitNodeReducedEq3}
\end{equation}
for $a_i\in K$.
\end{prp}

\begin{prp}[\cite{KH20b}, Proposition 3.3.1]\label{ReductionCusp}
Any trigonal curve of genus $5$ over $K$ in {\bf (Cusp case)} has a quintic model in $\mathbb{P}^2$ of the form \eqref{CuspReducedEq}:
%\begin{enumerate}
%\item[\rm (1)] 
\begin{equation}
\begin{split}
F &= x^2 z^3 + a_1 y^3 z^2 + 
(a_2 x^4 + a_3 x^3 y + a_4 x^2 y^2 + b_1 x y^3 + a_5 y^4) z \\
&\quad + a_6 x^5 + a_7 x^4 y + a_8 x^3 y^2 + a_9 x^2 y^3 + b_2 x y^4 + a_{10} y^5 \label{CuspReducedEq}
\end{split}
\end{equation}
%\end{enumerate}
for $a_i\in K$ $(i=1,\ldots, 10)$ with $a_1\ne 0$,
where $b_1 \in \{0,1\}$ and $b_2 \in \{0,1\}$.
\end{prp}

Here we recall a criterion on the superspecialty of the trigonal curve $C$ of genus $5$ in Proposition \ref{cor:genus5} below.
For the proof, see \cite[Section 2.2]{KH20b}, where the authors computed the Hasse-Witt matrix of $C$.

\begin{prp}[\cite{KH20b}, Corollary 2.2.2]\label{cor:genus5}
With notation as above, $C$ is superspecial if and only if the coefficients of the monomials $x^{pi-i'}y^{pj-j'}z^{pk-k'}$ in $F^{p-1}$ are equal to zero, where $(i,j,k)$ and $(i',j',k')$ run through $(3,1,1),(1,3,1),(2,2,1),(2,1,2),(1,2,2)$.
\end{prp}

Here we also describe a method given in \cite[Section 5.1]{KH20b} to test whether two trigonal curves of genus $5$ are $k$-isomorphic or not, where $k=K$ or $\overline{K}$.
Let $C_1$ and $C_2$ be trigonal curves of genus $5$ over $K$, and let $V(F_1)$ and $V(F_2)$ be the associate quintics in $\mathbb{P}^2$.
Recall from \cite[Lemma 2.1.2]{KH20b} that $C_1 \cong C_2$ over $k$ is equivalent to $V(F_1) \cong V (F_2)$ over $k$, i.e., there exist $M \in \mathrm{GL}_3 ( k )$ and $\lambda \in k^{\times}$ such that $M \cdot F_1 = \lambda F_2$.
Regarding entries of $M$ and $\lambda$ as variables, we reduce the (non-)existence of such $M$ and $\lambda$ into that of a solution over $k$ of a multivariate system, see \cite[Section 5.1]{KH20b} for more details including concrete algorithms to test the (non)-existence of such a solution.

Combining Propositions \ref{ReductionSplitNode} -- \ref{cor:genus5} with the isomorphism test described as above, we can construct a strategy similar to Strategy \ref{strategy1} for enumerating superspecial trigonal curves of genus $5$ over $K$, see \cite[Section 4]{KH20b} for concrete algorithms.
In Theorem \ref{MainTheorem2} below, we collect main results in \cite{KH20b} obtained by executing the algorithms over Magma.

\begin{thm}\label{MainTheorem2}
\begin{enumerate}
	\item (\cite[Theorem A]{KH20b}) There is no superspecial trigonal curve of genus $5$ in characteristic $7$.
	\item (\cite[Theorem B]{KH20b}) Any superspecial trigonal curve of genus $5$ over $\mathbb{F}_{11}$ is $\mathbb{F}_{11}$-isomorphic to the normalization of
\begin{equation}
x y z^3 + a_1 x^5 + a_2 y^5  = 0 \label{sscurve_F11}
\end{equation}
in $\mathbb{P}^2$, where $a_1, a_2 \in \mathbb{F}_{11}^{\times}$, or the normalization of
\begin{equation}
 (x^2 - 2 y^2) z^3 + a x^5 + b x^4 y + (9 a) x^3 y^2 + 4 b x^2 y^3 + ( 9 a ) x y^4 + 3 b y^5 =0 \label{sscurve_F11_2}
\end{equation}
in $\mathbb{P}^2$, where $(a, b) \in (\mathbb{F}_{11})^{\oplus 2} \smallsetminus \{ (0,0 )\}$.
	\item (\cite[Proposition 5.1.1 (I)]{KH20b}) There exist precisely four $\mathbb{F}_{11}$-isomorphism classes of superspecial trigonal curves of genus $5$ over $\mathbb{F}_{11}$.
Representatives of the four isomorphism classes are given by the normalization $C_i$ of $C_i^{\prime}=V (F_i) \subset \mathbb{P}^2$, where
\begin{eqnarray}
F_1 & = & x y z^3 + x^5 + y^5 , \nonumber \\
F_2 & = & x y z^3 + 2 x^5 + y^5 , \nonumber \\
F_3 & = & x y z^3 + 3 x^5 + y^5 , \nonumber \\
F_4 & = & (x^2 - 2 y^2) z^3  + x^5 + 9 x^3 y^2 + 9 x y^4. \nonumber
\end{eqnarray}
\item (\cite[Proposition 5.1.1 (II)]{KH20b}) There exists a unique $\overline{\mathbb{F}_{11}}$-isomorphism class of superspecial trigonal curves of genus $5$ over $\mathbb{F}_{11}$.
A representative of the unique isomorphism class is given by the normalization $C^{\rm (alc)}$ of the singular curve $C^{\rm (alc)^{\prime}} = V (F) \subset \mathbb{P}^2$ with $F =  x y z^3 + x^5 + y^5$.
	\item (\cite[Theorem C]{KH20b}) There is no superspecial trigonal curve of genus $5$ over $\mathbb{F}_{13}$.
\end{enumerate}
\end{thm}

In Table \ref{tab:genus5}, we summarize known values of the number $\# \Lambda_{5,p}$ of $\overline{\mathbb{F}_p}$-isomorphism classes of superspecial curves of genus $5$ in characteristic $p$.
As described at the beginning of this section, there is no superspecial non-hyperelliptic (resp.\ hyperelliptic) curve of genus $5$ for $p \leq 5$ (resp.\ $p \leq 7$).
The number written in bold type is determined by our theorem (Theorems \ref{MainTheorem2}) described in this section.
The notation `{\bf H}', `{\bf T}' and `{\bf G}' denote the hyperelliptic, trigonal, generic (canonical and non-trigonal) cases respectively.
The number written in each bracket is the number of $\overline{\mathbb{F}_p}$-isomorphism classes of  superspecial curves $C$ such that $C$ has a model over $\mathbb{F}_p$.

\begin{table}[h]
\centering
    \caption{Known values of the number $\# \Lambda_{5,p}$ of $\overline{\mathbb{F}_p}$-isomorphism classes of superspecial curves of genus $5$ in characteristic $p$.
}\label{tab:genus5}
  \begin{tabular}{|c||c|c|c|c|c|c|c|c|c|} \hline
  $p$                    & $2$   & $3$   & $5$  & $7$   & $11$     &  $13$ & $17$  & $19$ & $\geq 23$\\ \hline
  \multirow{2}{*}{\bf H} & \multirow{2}{*}{0} & \multirow{2}{*}{0} & \multirow{2}{*}{0} & \multirow{2}{*}{0} & \multirow{2}{*}{?} & \multirow{2}{*}{?} & \multirow{2}{*}{?} & \multirow{2}{*}{?} & \multirow{2}{*}{?} \\ 
                         &  & & & & & & & &  \\ \hline 
  \multirow{2}{*}{\bf T} & \multirow{2}{*}{0}   & \multirow{2}{*}{0} &  \multirow{2}{*}{0}  & \multirow{2}{*}{\bf 0}   & ? &  ? & \multirow{2}{*}{?} & \multirow{2}{*}{?}    & \multirow{2}{*}{?}\\ 
                         &  & &  & &({\bf 1})     & ({\bf 0}) & & &  \\ \hline 
\multirow{2}{*}{\bf G} & \multirow{2}{*}{0}   & \multirow{2}{*}{0} &  \multirow{2}{*}{0}  & \multirow{2}{*}{?}   & \multirow{2}{*}{?} &  \multirow{2}{*}{?} & \multirow{2}{*}{?} & \multirow{2}{*}{?}    & \multirow{2}{*}{?}\\ 
                         &  & &  & &     &      & & &  \\ \hline 
\multirow{2}{*}{$\# \Lambda_{5,p}$} & \multirow{2}{*}{0}   & \multirow{2}{*}{0} &  \multirow{2}{*}{0}  & \multirow{2}{*}{?}   & \multirow{2}{*}{?} & \multirow{2}{*}{?} & \multirow{2}{*}{?} & \multirow{2}{*}{?}    & \multirow{2}{*}{?}\\ 
                         &  & &  & &     &  & & &  \\ \hline 
  \end{tabular}
\end{table}

We close this section with open problems in the enumeration of superspecial curves of genus five:

\begin{prob}[Genus five]
Determine the number of $K$ or $\overline{K}$-isomorphism classes of superspecial curves of genus five over $K$ in the following cases:
\begin{enumerate}
\item Hyperelliptic case over $K = \mathbb{F}_{p}$ or over $K =\mathbb{F}_{p^2}$ for $p \geq 11$.
\item Trigonal case over $K = \mathbb{F}_{p}$ for $p \geq 17$ or over $K =\mathbb{F}_{p^2}$ for $p \geq 11$.
\item Generic (canonical and non-trigonal) case over $K = \mathbb{F}_{p}$ or over $K =\mathbb{F}_{p^2}$ for $p \geq 7$.
\end{enumerate}
\end{prob}

% \appendix

%=============================================================================
\section{Enumeration of certain genus-four superspecial curves}\label{sec:app}
%=============================================================================

While our papers \cite{KH17} and \cite{KH20} (resp.\ \cite{KH18}) enumerate superspecial curves among the {\it whole} space of non-hyperelliptic (resp.\ hyperelliptic) curves of genus $4$, the paper~\cite{KHH} enumerates those among a certain family of non-hyperelliptic curves of genus $4$, that is, {\it Howe curves}.
In \cite{Howe}, these curves were first studied in order to quickly construct genus-$4$ curves with many rational points, and also they heuristically tend to be superspecial.
It was also proved in \cite{KHS20} that there exists a supersingular Howe curve in every characteristic $p>3$.

In this section, we briefly describe results of \cite{KHH}, where the authors present computational methods ((A), (B) and (C) below) for enumerating superspecial Howe curves.
We start with recalling the definition of a Howe curve.
Throughout this section, let $K$ be an algebraically closed field of characteristic $p \neq 2$.

\begin{dfn}
A {\it Howe curve} over $K$ is a curve which is isomorphic to the normalization of the fiber product $E_1 \times_{\mathbb{P}^1} E_2$ of two genus-$1$ double covers $E_i \rightarrow {\mathbb{P}}^1$ ramified over $S_i$, where each $S_i$ consists of $4$ points and where $\# (S_1\cap S_2)=1$.
\end{dfn}

The superspeciality of a Howe curve is reduced into that of curves of low genera as follows:
For a Howe curve $H$ with two genus-$1$ double covers $E_i : y^2 = f_i(x) $, where $f_i$ is a separable polynomial of degree $3$ or $4$ with $i=1,2$, we have a genus $2$-curve $C : y^2 = f_1 f_2$.
It follows from \cite[Theorem 2.1]{Howe} (see also \cite[Theorem C]{KL89} for a more general result) that $H$ is supersingular if and only if $E_1$, $E_2$ and $C$ are all supersingular.

In order to enumerate superspecial Howe curves, two strategies (A) and (B) below are provided in \cite{KHH}.
The authors of \cite{KHH} also gave a method ((C) below) to decide whether two Howe curves are $K$-isomorphic, or not.

%==================================================================
\subsection*{(A) $(E_1,E_2)$-first, using Cartier-Manin matrices}
%==================================================================
In this strategy, we use the same realization of Howe curves as in \cite{KHS20}, that is, the fiber product of $E_1\colon z^2 = f_1(x):=x^3 + A_1 \mu^2 x + B_1 \mu^3$ and $E_2\colon w^2 = f_2(x):= (x-\lambda)^3 + A_2 \mu^2 (x-\lambda) + B_2 \mu^3$ over $\mathbb{P}^1=\mathrm{Proj} (K[x,y])$, where $A_1$, $B_1$, $A_2$ and $B_2$ are elements in $K$ such that $E_{A_i,B_i} : y^2 = x^3 + A_i x + B_i$ ($i=1,2$) are supersingular elliptic curves over $K$, and where $\lambda$, $\mu$ and $\nu$ are elements in $K$ such that (i) $\mu \neq 0$ and $\nu \neq 0$, and (ii) $f_1$ and $f_2$ are coprime. 
Note that a point $(\lambda : \mu : \nu) \in \mathbb{P}^2(K)$ satisfying (i) and (ii) is said to be {\it of Howe type} in \cite{KHS20}.
It was shown in \cite[Proposition 4.1]{KHH} that any superspecial Howe curve is $K$-isomorphic to the normalization of $E_1 \times_{\mathbb{P}^1} E_2$ obtained as above for $A_1$, $B_1$, $A_2$, $B_2$, $\lambda$, $\mu$ and $\nu$ belonging to $\mathbb{F}_{p^2}$.

This strategy enumerates pairs of supersingular elliptic curves $E_i : y^2 = f_i (x)$ ($i=1,2$) so that $C : y^2 = f_1(x) f_2(x)$ is superspecial.
To do this, for each unordered pair of $(A_1,B_1)$ and $(A_2,B_2)$, it suffices to compute the solutions $(\lambda : \mu : \nu)\in \mathbb{P}^2(\mathbb{F}_{p^2})$ (of Howe type) to the homogeneous system $M \equiv 0$, where $M$ is the Cartier-Manin matrix of $C$.
Once all pairs $(E_1,E_2)$ are enumerated, we classify isomorphism classes of Howe curves defined by the pairs, by the isomorphism test described in (C) below.
For a concrete algorithm of the enumeration based on this strategy, see \cite[Section 4.2]{KHH}.
It is also shown in \cite[Section 4.3]{KHH} that the complexity of this algorithm is $\tilde{O}(p^6)$ arithmetic operations in $\mathbb{F}_{p^2}$.

%====================================================
\subsection*{(B) $C$-first, using Richelot isogenies}
%====================================================
The second strategy first enumerates superspecial curves $C\colon y^2=f(x)$ of genus $2$, where $f(x) \in K[x]$ is a separable polynomial of degree $6$.
Once all superspecial curves of genus $2$ are enumerated, we then enumerate decompositions $f(x)=f_1(x)f_2(x)$ with $f_i(x)$ of degree $3$ so that there is an element $b \in K \cup \{ \infty \}$ that makes both genus-$1$ curves $E_i\colon y^2 = (x-b)f_i(x)$ ($i=1,2$) supersingular.
% The moduli space of curves of genus $2$ is of dimension $3$.
% Since this dimension is larger than the space of $(\lambda:\mu:\nu)\in \mathbb{P}^2$ considered in (A), this strategy, a priori, looks inefficient.
% However, surprisingly, we conclude that strategy (B) enumerates superspecial Howe curves much more efficiently than does (A).
% The advantage of (B) comes from making use of Richelot isogenies.

For this, we first apply a method given in \cite[Section 3]{HLP} to construct some (at least one) superspecial curves of genus $2$ by gluing supersingular elliptic curves together along their $2$-torsion.
We then produce more such curves by applying {\it Richelot isogenies} to the curves already produced, where the definition of a Richelot isogeny of two curves of genus $2$ is as follows:

\begin{dfn}
Two genus-$2$ curves are {\it Richelot isogenous} if there exits an isogeny $\Psi : J(C_1) \rightarrow J(C_2)$ such that $\mathrm{Ker}(\Psi)$ is isomorphic to $\mathbb{Z}/2\mathbb{Z} \times \mathbb{Z}/2\mathbb{Z}$ that is maximal isotropic with respect to $2$-Weil pairing.
In this case, the isogeny $\Psi$ is called a {\it Richelot isogeny}. 
\end{dfn}

Note that given a curve $C_1$, we can compute (at most $15$) genus-$2$ curves $C_2$, which are Richelot isogenous to $C_1$, with an isogeny map $\Psi$, see e.g., \cite[Chapter 8]{Smith} for more details.
The above procedure to enumerate superspecial genus-$2$ curves terminates because there are only finitely many superspecial curves of genus $2$, and a recent result of Jordan and Zaytman~\cite[Theorem 43]{JZ} shows that we obtain all isomorphism classes of superspecial curves of genus $2$ in this way.
Once all pairs $(E_1,E_2)$ are enumerated, we classify isomorphism classes of Howe curves defined by the pairs, by the isomorphism test described in (C) below.
A concrete algorithm of the enumeration of superspecial Howe curves based on this strategy is given in \cite[Section 5.3]{KHH}, and its complexity is $\tilde{O}(p^4)$ arithmetic operations in $\mathbb{F}_{p^2}$.

%========================================================
\subsection*{(C) A new isomorphism test for Howe curves.}
%========================================================
% Strategy (A) above produces many not-necessarily-distinct Howe curves, so to prevent overcounting we are left with the task ofproducing a unique representative for each isomorphism class.
Since every Howe curve is canonical (cf.\ \cite[Lemma 2.1]{KHH}), we can test whether two Howe curves are isomorphic or not, by applying the isomorphism test for canonical curves of genus $4$ given in Section \ref{sec:genus4} (cf.\ \cite[Section 6.1]{KH17} and \cite[Section 4.3]{KH18}).
However, this turns out to be very costly since it uses many Gr\"obner basis computations.
In \cite[Section 3]{KHH}, the authors present an efficient isomorphism test specific to Howe curves.
We here briefly describe the isomorphism test given in \cite[Section 3]{KHH}.

We first recall from the beginning of \cite[Section 3]{KHH} that a Howe curve is specified by the following three pieces of information:
(1) A genus-$2$ curve $C$.
(2) An unordered pair of disjoint sets $\{W_1, W_2\}$, each consisting of three Weierstrass points of $C$.
(3) An unordered pair of distinct points $\{P_1, P_2\}$ on $C$ that are mapped to one another by the hyperelliptic involution.
We here call $(C, \{W_1, W_2\}, \{P_1, P_2\})$ {\it a Howe triple} of a Howe curve.
A criterion given in \cite[Section 3]{KHH} for determining whether two Howe curves are isomorphic or not is the following:

\begin{prp}[\cite{KHH}, Corollary 3.3]\label{C:isomorphism}
Two Howe triples $(C, \{W_1, W_2\}, \{P_1, P_2\})$ and $(C', \{W_1', W_2'\}, \{P_1', P_2'\})$ give isomorphic Howe curves if and only if there is an isomorphism $C \to C'$ that takes $\{W_1, W_2\}$ to $\{W_1', W_2'\}$ and $\{P_1, P_2\}$ to $\{P_1', P_2'\}$. 
\end{prp}

This isomorphism test is conducted by simply deciding whether there exist any automorphisms of $\mathbb{P}^1$ that respect the sets of Weierstrass points and their divisions, and that take the $x$-coordinate of $P_1$ and $P_2$ to that of $P_1'$ and $P_2'$.
Clearly this procedure does not require any Gr\"obner basis computation, and also it is shown to be more efficient than the isomorphism test for canonical curves of genus $4$ given in Section \ref{sec:genus4}.

%===========================================================
\subsection*{Main theorems in \cite{KHH} and open questions}
%===========================================================

The authors of \cite{KHH} implemented algorithms based on (A) -- (C) over Magma, and executed them to enumerate superspecial Howe curves for concrete $p$.
Recall that the complexities of (A) and (B) are $\tilde{O}(p^6)$ and $\tilde{O}(p^4)$ respectively.  
Practical time behaver of (A) and (B) for $5 \leq p \leq 53$ is shown in \cite[Table 2]{KHH}.
As the estimated complexities show, we expect from \cite[Table 2]{KHH} that (B) is extremely faster than (A) in practice; e.g., for $p=53$, (A) takes $5678.32$ seconds, while (B) takes only $1.46$ seconds, under the authors' experimental environment (details are written in \cite[Section 6]{KHH}).
From this, the authors of \cite{KHH} decided to adopt (B) in order to obtain results for $p$ larger than $53$.

Main results obtained by the execution of (B) together with (C) are the following:

\begin{thm}
\begin{enumerate}
\item (\cite[Theorem 1.1]{KHH}) For every prime $p$ with $7 < p < 20000$, there exists a superspecial Howe curve in characteristic $p$.
\item (\cite[Theorem 1.2]{KHH}) For every prime $p$ with $7 < p \le 199$, the number of $\overline{\mathbb{F}_{p}}$-isomorphism classes of superspecial Howe curves in characteristic $p$ is given in Table \textup{\ref{table:enumerate}}.
\end{enumerate}
\end{thm}

\begin{table}[ht]
\caption{For each prime $p$ from $11$ to $199$, we give the number $n(p)$ of superspecial Howe curves over $\overline{\mathbb{F}_p}$ and the ratio of $n(p)$ 
to the heuristic prediction $p^3/1152$ (see \cite[Section 5]{KHH}).}
\label{table:enumerate}
\vspace{-1ex}
\begin{center}
\renewcommand{\arraystretch}{0.9}
\begin{tabular}{|rrrrrrrrrrr|}
% \toprule
\hline
  $p$ &  $n(p)$ &\small    Ratio &\hbox to 1em{}&   $p$ &  $n(p)$ &    Ratio &\hbox to 1em{}&   $p$ &  $n(p)$ &    Ratio \\ \hline
% \cmidrule{1-3}\cmidrule{5-7}\cmidrule{9-11}
 $11$ &     $4$ &  $3.462$ &&  $67$ &   $260$ &  $0.996$ && $137$ &  $2430$ &  $1.089$ \\ 
 $13$ &     $3$ &  $1.573$ &&  $71$ &   $742$ &  $2.388$ && $139$ &  $2447$ &  $1.050$ \\
 $17$ &    $10$ &  $2.345$ &&  $73$ &   $316$ &  $0.936$ && $149$ &  $3082$ &  $1.073$ \\
 $19$ &     $4$ &  $0.672$ &&  $79$ &   $595$ &  $1.390$ && $151$ &  $3553$ &  $1.189$ \\
 $23$ &    $33$ &  $3.125$ &&  $83$ &   $655$ &  $1.320$ && $157$ &  $3427$ &  $1.020$ \\
 $29$ &    $45$ &  $2.126$ &&  $89$ &   $863$ &  $1.410$ && $163$ &  $3518$ &  $0.936$ \\
 $31$ &    $59$ &  $2.281$ &&  $97$ &   $802$ &  $1.012$ && $167$ &  $6268$ &  $1.550$ \\
 $37$ &    $41$ &  $0.932$ && $101$ &  $1207$ &  $1.350$ && $173$ &  $4780$ &  $1.064$ \\
 $41$ &   $105$ &  $1.755$ && $103$ &  $1151$ &  $1.213$ && $179$ &  $5771$ &  $1.159$ \\
 $43$ &    $79$ &  $1.145$ && $107$ &  $1237$ &  $1.163$ && $181$ &  $5419$ &  $1.053$ \\
 $47$ &   $235$ &  $2.608$ && $109$ &  $1193$ &  $1.061$ && $191$ &  $9610$ &  $1.589$ \\
 $53$ &   $167$ &  $1.292$ && $113$ &  $1323$ &  $1.056$ && $193$ &  $6298$ &  $1.009$ \\
 $59$ &   $259$ &  $1.453$ && $127$ &  $2013$ &  $1.132$ && $197$ &  $6839$ &  $1.030$ \\
 $61$ &   $243$ &  $1.233$ && $131$ &  $2606$ &  $1.335$ && $199$ &  $8351$ &  $1.221$ \\ \hline
% \bottomrule
\end{tabular}
% \vskip 3ex
\end{center}
\end{table}

We can easily increase the upper bounds on $p$ in these two theorems.
For example, on a 2.8 GHz Quad-Core Intel Core i7 with 16GB RAM, computing the $8351$ superspecial Howe curves in characteristic $199$ using an algorithm based on (B) took $124$ seconds in Magma.
Finding examples of superspecial Howe curves for every $p$ between $7$ and $20000$ took $680$ minutes on the same PC.

We close this section with an open problem in the enumeration of superspecial Howe curves, and an open question on the existence of such curves:

\begin{prob}
\begin{enumerate}
\item (Problem) Determine the number of $\mathbb{F}_{p^2}$-isomorphism classes of superspecial Howe curves.
\item (Question) Does there exist a superspecial Howe curve in any characteristic $p > 7$?
\end{enumerate}
\end{prob}

\noindent Department of Mathematical Informatics, Graduate School of Information Science and Technology, The University of Tokyo, Hongo 7-3-1, Bunkyo-ku, Tokyo, 113-8656, Japan.

\noindent E-mail: \texttt{kudo@mist.i.u-tokyo.ac.jp}
% \end{flushright}

\end{document}